\documentclass[12pt]{article}
\usepackage{fullpage}
\usepackage{epsf}
\newcommand {\eq}[1]{\begin{equation}\label{#1}}
\newcommand {\en} {\end{equation}}
\newcommand {\proof} {\par{\it Proof}. \ignorespaces}
\newcommand {\eproof}
      {\space
        {\ \vbox{\hrule\hbox{\vrule height1.3ex\hskip0.8ex\vrule}\hrule}}
        \par}
\newcommand {\mat}[1]{\left[\begin{array}{#1}}
\newcommand {\rix}          {\end{array}\right]}

\newtheorem{theorem}          {Theorem}

\newtheorem{corollary}     [theorem]{Corollary}
\newtheorem{proposition}     [theorem]{Proposition}

\newtheorem{example}        [theorem]   {Example}

\newcommand {\diag}     {\mathop{\rm diag}\nolimits}

\newcommand{\lam}       {\{ \lambda_{1},\ldots,\lambda_{n}\}}

\def\eqbd{\mathop{{:}{=}}}

\def\C{{\rm C\kern-.48em\vrule width.06em height.6em depth-.02em
                 \kern.48em}}
\def\om{\omega}
\def\lam{\lambda}

\def\Zbar{\overline{Z}}

\date{13.04.2003}
\title{Potter, Wielandt, and Drazin
on the matrix equation $AB = \omega BA$,
with some new answers to old questions}

\author{Olga Holtz \thanks{
         Inst. f. Mathematik, MA 4-5,
         Technische Universit\"at Berlin,
         D-10623 Berlin,
         Fed. Rep. Germany. E-mail: holtz@math.TU-Berlin.DE\/
On leave from CS Dept., Univ. of Wisconsin, Madison, WI 53706, USA.
Supported by Alexander von Humboldt Foundation.}\and
        Volker Mehrmann \thanks{
         Inst. f. Mathematik, MA 4-5,
         Technische Universit\"at Berlin,
         D-10623 Berlin,
         Fed. Rep. Germany. E-mail: mehrmann@math.TU-Berlin.DE\/
Supported by DFG research grant Me790/15.}
                 \and
        Hans Schneider \thanks{ Dept. of Mathematics, Univ. of
Wisconsin, Madison, WI 53706, USA.
         E-mail: hans@math.wisc.edu \/Supported by DFG research grant Me790/15.
}
}
\begin{document}
\maketitle

\begin{abstract} In this partly historical and partly research oriented
note, we display a page of an unpublished mathematical diary of Helmut
Wielandt's for 1951. There he gives a new proof of a theorem due
to H.~S.~A.~Potter on the matrix equation $AB = \omega BA$, which is
related to the $q$-binomial theorem, and asks some further questions,
which we answer. We also describe results by M.~P.~Drazin and others
on this equation.
\end{abstract}

\section{Introduction}

The aim of this paper is to present a slice of the linear algebra
of the 1950's and to give some answers to questions raised then.

It was Helmut Wielandt's habit over many years to make notes in
what he called diaries (Tageb\"ucher) on papers that interested
him. Many notes are essentially summaries of a paper, but in other
cases Wielandt would add questions, ideas, or even further
results. In this note we discuss one such entry which appears on
page 35 of Diary VII (1951) which will appear in transcribed
electronic form \cite{WieVII51}. The entry concerns a paper which
Wielandt reviewed for the Zentralblatt. We next turn to this
paper.

In 1950, H. S. A.~Potter, a mathematician at Aberdeen University
in Scotland, published a note in the American Mathematical
Monthly \cite{Pot50}, on the matrix equation
  \eq{quasicomm}  AB = \omega BA. \en
He called a pair of  complex $n \times n$ matrices $A, B$
satisfying (\ref{quasicomm}) {\em quasi-commutative}.
We shall call matrices (\ref{quasicomm}) $\omega$-commutative, see Section
\ref{potpro} for a definition of this term applicable to general
rings. Otherwise we follow Potter's notation. It should be noted
here that the term "quasicommutative" has also been used in a
different sense, see~\cite{McC34}.

\begin{figure}
\epsfysize=4.3cm \epsfbox{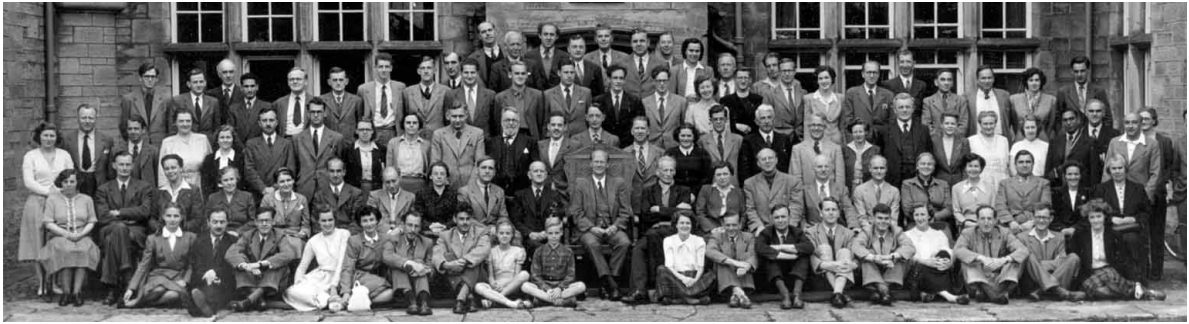}

\caption{\hskip 0.4cm H.~S.~A.~Potter is third from the right in the front row.
A.~C.~Aitken and H.~W.~Turnbull are tenth and eleventh from the left,
respectively, in the second row (seated).  H.~Schneider is the fourteenth
from the left in the fourth row.
 \label{figure1}}
\end{figure}


Potter's principal result is the following theorem:
\begin{theorem}[Potter~\cite{Pot50}] \label{pot}
Let $A$ and $B$ be complex square matrices satisfying
(\ref{quasicomm}) where $\omega$ is a primitive $q$-th root of
unity. Then \eq{ident} A^q + B^q = (A+B)^q. \en
\end{theorem}

In his note, Potter proves his theorem by deriving it from the
general expansion of $(x+y)^q$ for any nonnegative integer $q$ and
$\omega$-commutative $x$ and $y$ for arbitrary complex $\om$. This
formula, which we state as (\ref{genident})--(\ref{qfacs}),
involves the classical $q$-binomial coefficients and is
currently referred to as the noncommutative $q$-binomial theorem,
see, e.g.,~\cite[Formula 10.0.2]{AndAR99}
or~\cite[Exercise 1.35]{GasR90}. [But care: the
$q$ in the last sentence is our $\om$.]. The result that
(\ref{genident})--(\ref{qfacs}) holds for $\omega$-commutative operators is generally
attributed to Sch\"utzenberger \cite{Sch53}. We shall call
(\ref{genident})--(\ref{qfacs}) the {\em Potter-Sch\"utzenberger} formula.
It is
of considerable interest in the study of quantum groups, see for
example \cite[p.75]{Kas95}. In fact, Potter's proof shows that it
holds under very general conditions, which we examine in
Section~\ref{potpro}.

Potter refers and applies results in the book by
Turnbull-Aitken~\cite[p.148]{TurA32} where a matrix $X$
satisfying $AX = XC$ is called  {\em commutant} of $A$ and $C$.
There all commutants of $A$ and $C$ are determined on the
assumption that $A$ and $C$ are in Jordan canonical form. If $A$
and $B$ are quasi-commutative, then clearly $B$ is a commutant of
$A$ and $\omega A$. The general question of commutants was also
considered by Goddard-Schneider \cite{GodS55}. One might observe
that all the mathematicians mentioned in this paragraph were in
Scotland in the early 1950's. Figure~\ref{figure1} shows
the participants of the  1951 Edinburgh Mathematical Society Colloquium
at St.Andrews \cite{photo}. There four mathematicians mentioned in our
article are present.

Wielandt's proof of Potter's Theorem \ref{pot} is reproduced and
translated in Section~\ref{wnote}. We comment on it and give a
variant in Section~\ref{anotherproof}. This proof uses matrix
theory non-trivially and it is based on an insightful observation.
However, it heavily uses the assumption that $\om$ is a primitive
$q$-th root of $1$ and there is no obvious way of obtaining the
more general Theorem \ref{genpotter} using his methods.

In his diary, following the proof of Potter's theorem, Wielandt
also raises some questions. These include the construction of all
identities satisfied by $\omega$-commutative  matrices and the determination of
all irreducible pairs of $\omega$-commutative matrices. Naturally
unaware of Wielandt's question, M. P.~Drazin, then at Cambridge,
England, essentially answers the latter question in \cite{Dra51}.

In Section~\ref{nf} we take up the question of normal forms for
pairs of quasi-commutative matrices. We present the pre-normal
form obtained by Drazin~\cite{Dra51} and show that the
classification problem of quasi-commutative matrices is equivalent
to the classification problem of pairs of commuting matrices, both
under simultaneous similarity.

In Section~\ref{potterconv} we present counterexamples showing
that the converse to Potter's theorem does not hold, not even
for some of its weakened versions.

In Section~\ref{identities} we determine all polynomial
identities satisfied by quasi-commutative matrices thus answering
Wielandt's first question.

Finally, in section~\ref{final} we discuss work on quasi-commutative
matrices preceding that of Potter and Wielandt.

\section{Potter's proof}\label{potpro}
We begin by examining Potter's proof of Theorem \ref{pot}. In the
first part of the proof Potter does not assume that $\om$ is a
root of unity and for $\omega$-commutative
 matrices $A,B$ he proves the general formula
(here stated in a slightly different but equivalent form)
 \eq{genident} (A+B)^q = \sum_{k=0}^q c_k B^kA^{q-k} \en
where the $c_k$ are determined by
 \eq{binoms} \phi_k \phi_{q-k} c_k = \phi_q \quad, k = 0,\ldots,q, \en
 and the $\phi_k$ are given by
 \eq{qfacs} \phi_k = \prod_{s=1}^k (1+\cdots+\om^{s-1}) \quad, k =
 0,\ldots,q.
\en

 The coefficients $c_k$  in (\ref{binoms}), known as the
\emph{ $q$-binomial coefficients}, were well-studied in the
nineteenth century in the theory of hypergeometric series, see for
example \cite[Chapter 10]{AndAR99} and in the theory of partitions
combinatorics, see~\cite[Chapter 11]{AndAR99} and \cite[Section
1.3]{Sta86}.

Let $R$ be any ring with identity $1$ and let $\om,x$ and $y$ be
elements of $R$. Let $\Zbar$ be the subring generated by $1$ in
$R$. Thus $\Zbar$ is isomorphic either to the ring of integers $Z$
or $Z_m$ the ring of integers mod $m \in Z$. We call $x$ and
$y$ {\em $\omega$-commutative} if the following identities hold
 \begin{eqnarray}
 \om x &=& x \om \nonumber \\
 \om y &=& y \om \nonumber \\
 xy &=& \om yx.
\label{basicids}
 \end{eqnarray}
 By Potter's argument we may obtain the following version of the
Potter-Sch\"utzenberg theorem.
\begin{theorem} \label{genpotter}
 Let $R$ be a ring with $1$ and let $\Zbar$ be the subring
 generated by $1$. Let $\om \in R$ and let $x$ and $y$ be $\omega$-commutative
 elements of $R$.
 Then
 \eq{genident2} (x+y)^q = \sum_{k=0}^q b_k y^kx^{q-k}, \en
where the $c_k$ and $\phi_k, \ k = 0,\ldots,q$, are given by
(\ref{binoms}) and (\ref{qfacs}) respectively.
\end{theorem}
 We observe that the coefficients $b_k$ lie in $\Zbar[\om]$ and
thus there is no loss of generality by considering only the
subring $\Zbar[\om,x,y]$ of $R$.
\begin{corollary} \label{speccase} Suppose that $\Zbar[\om]$ is
an integral
domain. Under the conditions of Theorem~\ref{genpotter}, suppose
further that
 \eq{nonzero} \phi_k \neq 0, \quad k = 1,\ldots, q-1, \en
 but that
\eq{zero} \phi_q = 0 . \en
 Then
 \eq{specid} (x+y)^q = x^q + y^q. \en
\end{corollary}
Evidently, if $\Zbar[\om]$ is a field and $\om$ is a primitive
$q$-th root of $1$ $R$, then (\ref{nonzero}) and (\ref{zero})
hold. These also hold if $\Zbar = Z_q$ and $\om = 1$. We may also
note that in the case of an integral domain $\Zbar[\om]$ a
necessary condition for (\ref{specid}) to be satisfied is that
(\ref{zero}) holds.

Let $K[x,y]$ be the ring in two noncommutative indeterminates $x$
and $y$ over a central field $K$. If $x$ and $y$ are  subject to
the relation $xy=\om yx$ where $\om\in K$
then $K[x,y]$ is today called a {\it quantum plane} over $K$, see
\cite[p.72]{Kas95}.

\section{Wielandt's notes}\label{wnote}
\begin{figure}
\epsfysize=20.cm \epsfbox{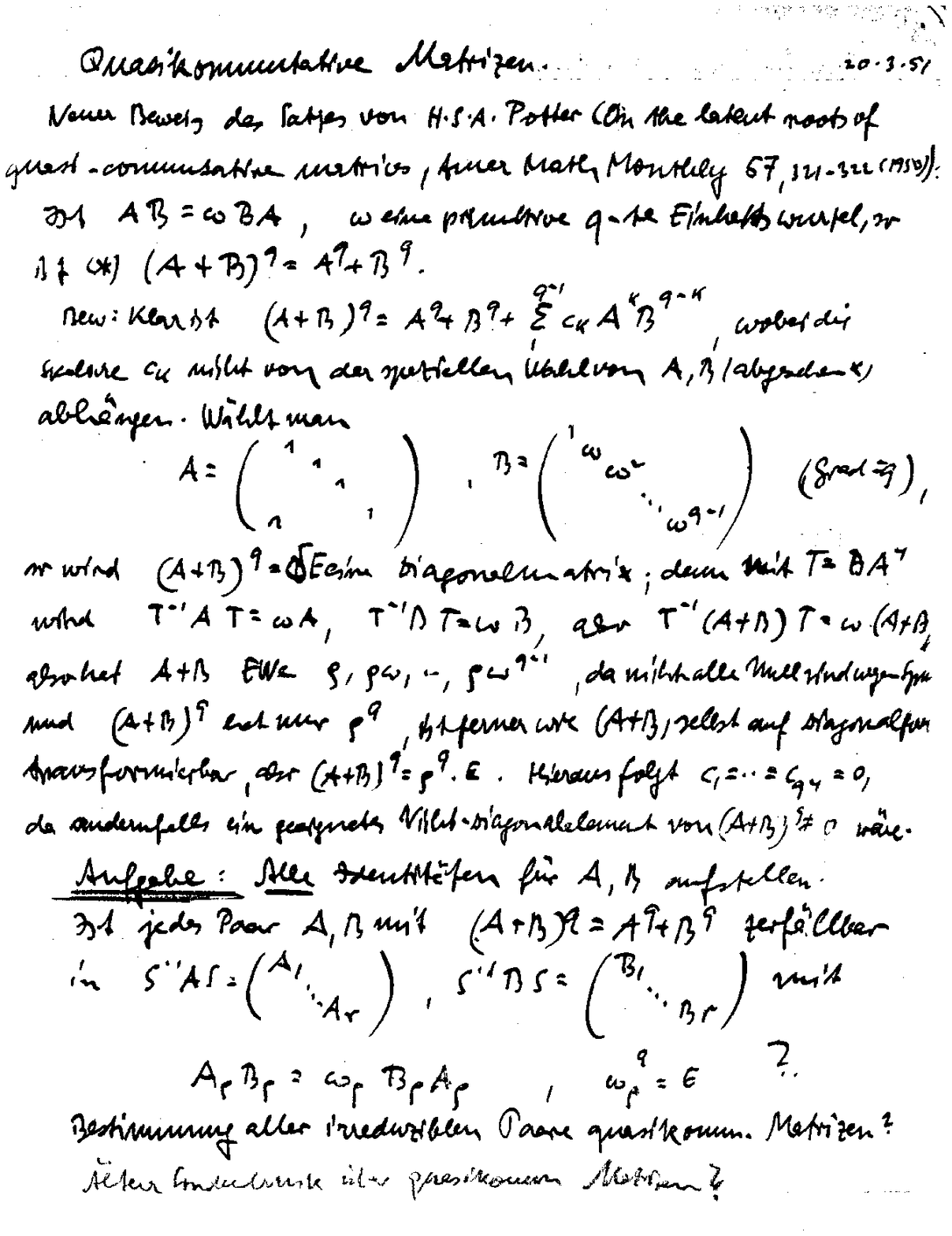}

\caption{\hskip 0.4cm Page 35 of Wielandt's diary VII.
\label{fax}}
\end{figure}

In Figure~\ref{fax} we display a facsimile of page 35 of
Wielandt's Diary VII, \cite{WieVII51}, dated 20 March 1951.

The transcription reads as follows:
{\it

Neuer Beweis des Satzes von H.~S.~A.~Potter  (On the latent roots of
quasi-commutative matrices, Amer. Math. Monthly 57, 321--322 (1950)).

Ist $AB= \omega BA$, $\omega$ eine primitive $q$-te Einheitswurzel, so
ist (*) $(A+B)^q=A^q+B^q$.

Bew: Klar ist $(A+B)^q=A^q+B^q+ \sum_1^{q-1} c_k A^kB^{q-k}$,
wobei die Skalare $c_k$ nicht von der speziellen Wahl von $A,B$
(abgesehen *) abh\"angen. W\"ahlt man

$$ A=\left ( \begin{array}{ccccc} & 1 &&& \\ && 1 && \\ &&& 1 & \\
&&&& 1 \\ 1&&&& \end{array} \right ),
\quad B= \left ( \begin{array}{ccccc} 1&  &&& \\ &\omega & && \\
    &&\omega^2&  & \\
&&&\ddots&  \\ &&&& \omega^{q-1} \end{array} \right ), (grad =q)
$$
 so wird  $(A+B)^q=\delta E$ eine Diagonalmatrix;
denn mit $T=BA^{-1}$ wird $T^{-1} A T =\omega A$, $T^{-1} B T =\omega
B$, also $T^{-1} (A+B) T =\omega (A+B)$, also hat
$A+B$ Ewe $\rho, \rho \omega, \ldots , \rho \omega^{q-1}$, da nicht
alle Null sind wegen Spur, und
$(A+B)^q$ hat nur $\rho^q$, ist ferner wie $(A+B)$ selbst auf
Diagonalform transformierbar, aber $(A+B)^q=\rho^q E$. Hieraus folgt
$c_1=\ldots=c_{q-1}=0$, da andernfalls  ein geeignetes
Nicht-Diagonalelement von $(A+B)^q\neq 0$ w\"are.

\underline {Aufgabe:} Alle Identit\"aten f\"ur $A,B$ aufstellen.
Ist jedes Paar $A,B$ mit $(A+B)^q=A^q+B^q$ zerf\"allbar in

$$ S^{-1} A S = \left ( \begin{array}{ccc} A_1 && \\
& \ddots & \\ && A_r \end{array} \right ), \quad
S^{-1} B S = \left ( \begin{array}{ccc} B_1 && \\
& \ddots & \\ && B_r \end{array} \right ) $$
mit
$$A_\rho B_\rho =\omega_\rho B_\rho A_\rho, \quad \omega_\rho^q=E
? $$
\\
Bestimmung aller irreduziblen Paare quasikommutativer Matrizen ?
\\
\"Altere Sonderdrucke \"uber quasikommutative Matrizen? }
\vskip .5truecm

The following is a translation of this note.
 \vskip .5truecm

{\it Quasi-commutative Matrices

New proof of a theorem of H.~S.~A.~Potter  (On the latent roots of
quasi-commutative matrices, Amer. Math. Monthly 57, 321--322 (1950)).

If $AB= \omega BA$, and $\omega$  is a primitive $q$-th root of unity,
then  (*) $(A+B)^q=A^q+B^q$.

Proof: It is clear that $(A+B)^q=A^q+B^q+ \sum_1^{q-1} c_k A^kB^{q-k}$,
where the  scalars $c_k$ do not depend on the special choice of of $A,B$
(except for *). If one chooses

$$ A=\left ( \begin{array}{ccccc} & 1 &&& \\ && 1 && \\ &&& 1 & \\
&&&& 1 \\ 1&&&& \end{array} \right ),
\quad B= \left ( \begin{array}{ccccc} 1&  &&& \\ &\omega & && \\
    &&\omega^2&  & \\
&&&\ddots&  \\ &&&& \omega^{q-1} \end{array} \right ), \quad (degree=q)
$$
then  $(A+B)^q=\delta E$ is a diagonal matrix;
since with $T=BA^{-1}$ also $T^{-1} A T =\omega A$, $T^{-1} B T =\omega
B$, and thus $T^{-1} (A+B) T =\omega (A+B)$, hence
$A+B$ has eigenvalues $\rho, \rho \omega, \ldots , \rho \omega^{q-1}$,
since not all are zero due to the trace and
$(A+B)^q$ has only $\rho^q$ as eigenvalue and is as $(A+B)$ transformable
to diagonal form, but $(A+B)^q=\rho^q E$. This implies that
$c_1=\ldots=c_{q-1}=0$, since otherwise a particular
non-diagonal element would satisfy  $(A+B)^q\neq 0$.

\underline {Problem:} Determine all identities for  $A,B$.
Is every pair  $A,B$ with $(A+B)^q=A^q+B^q$ decomposable as
$$
 S^{-1} A S = \left ( \begin{array}{ccc} A_1 && \\
& \ddots & \\ && A_r \end{array} \right ),  \quad
S^{-1} B S = \left ( \begin{array}{ccc} B_1 && \\
& \ddots & \\ && B_r \end{array} \right )
$$
with
$$A_\rho B_\rho =\omega_\rho B_\rho A_\rho, \quad \omega_\rho^q=E
? $$
\\
Determination of all irreducible pairs of quasi-commutative
matrices?
\\
Earlier work on quasi-commutative matrices? }
\vskip 0.5cm

Apparently, in the first paragraph of the proof Wielandt means
the quasi-commutativity relation rather than the relation~(*)
when he says `except for (*)'. Also, he must mean that $\omega_\rho^q=1$
since he is considering the case when $\omega_\rho$ is a scalar.


\section{Wielandt's proof and a variant}\label{anotherproof}

Wielandt's proof begins with the simple but insightful remark that
for $\omega$-commutative  matrices the coefficients $c_k$ in the expansion
$(A+B)^q = \sum_{k=0}^q c_k B^kA^{q-k}$ are independent of the
particular matrices $A,B$, and hence the result is proved if he
can show that the coefficients must be $0$ in the case of a
well-chosen pair of matrices $A$ and $B$.
The argument requires the linear independence of the set of
matrices $B^k A^{n-k},\ k = 1,\ldots,q-1$. Though Wielandt
does not say this, he chooses
a pair of matrices $A,B$ that satisfy this condition. He then
uses an argument involving eigenvalues  and the diagonability of
matrices to show that $c_1 = \ldots = c_{q-1} = 0$.

We now give a variant of Wielandt's proof.
Let $A$ and $B$ be the matrices chosen by Wielandt and let $s$ and
$t$ be any complex numbers. Since the eigenvalues of $B$ are the
$q$-th roots of unity, it follows that the characteristic
polynomial of $sB$ is $\lambda^q -s^q$. Since the proper principal
minors of $sA +tB$ and $tB$ coincide and $\det(sA+tB) =(-1)^{q-1}
(s^q+t^q)$, it follows that the characteristic polynomial of
$sA+tB$ is $\lambda^q - (s^q+t^q)$. By the Cayley-Hamilton theorem
\cite{Gan59} we obtain
\begin{equation}
(sA+tB)^q = (s^q+t^q)I = (sA)^q +(tB)^q.\label{identity}
\end{equation}
But $(sA+tB)^q=(sA)^q+(tB)^q+ \sum_{k=1}^{q-1} c_k s^kt^{q-k}A^kB^{q-k}, $
so each matrix coefficient $c_k A^kB^{q-k}$ must equal zero, which
implies $c_k=0$, since $A$ and $B$ are both nonsingular.
\eproof

Thus, the alternative proof demonstrates the following extension of
Potter's Theorem.
\begin{proposition}\label{potext}
Let $A$ and $B$ be quasi-commutative matrices satisfying
(\ref{quasicomm}) where $\omega$ is a primitive $q$-th root of
unity. Then \eq {fermat} (sA+tB)^q=(sA)^q+(tB)^q \en for all
$s,t\in \C$.
\end{proposition}

A proof in a rather similar spirit is given by R. Bhatia and L. Elsner
in~\cite{BhaE93} for the following fact: Let $A$ and $B$ be quasicommutative,
then the spectrum of $A+B$ is $p$-Carollian, i.e., the eigenvalues of
$A+B$ can be enumerated as
$$ (\lambda_1, \ldots, \lambda_r, \omega \lambda_1,\ldots, \omega \lambda_r,
\ldots, \omega^{p-1} \lambda_1,\ldots, \omega^{p-1} \lambda_r  ).$$
Moreover, the same holds for all perturbations of $B$ of specific
form given in~\cite[Theorem~2]{BhaE93}. The term `Carollian' was
invented by R.~Bhatia in honor of L. Carrol, initially to denote an
$n$-tuple that contains $-x$ if it contains $x$, and later turned
into `$p$-Carollian' for $n$-tuples
that contain all multiples of $x$  with $p$-th roots of unity. It is
used also in~\cite{BhaE92} and~\cite{Bha97}.

\section{Normal forms for quasi-commutative matrices}\label{nf}

Note that Wielandt asks the question on classification of
irreducible quasi-commutative pairs, having in mind reductions by
simultaneous similarity $A\mapsto T^{-1}AT$, $B\mapsto T^{-1}BT$,
which leave the relation~(\ref{quasicomm}) invariant.

To study this question, we start with some preliminary observations.
Suppose that
\eq{gqc}
AB=\alpha BA,
\en
where  $\alpha$ is a nonzero complex number.
By the above remark on simultaneous similarity,
we may assume w.l.o.g that  $A$ is in Jordan canonical form
$$A =\diag(J_1(\lam_1),\ldots,J_s(\lam_s)) $$
where $J_i(\lam_i)$ is a Jordan
block of size $k_i$ corresponding to the eigenvalue $\lam_i,\ i = 1,
\ldots,s$. We partition $B$ conformably with $A$, where  $B_{ij}$ is a
block of size $k_i \times k_j$. Using the construction
of \cite[p.148]{TurA32}, see also \cite{Gan59},
we conclude that  that $B_{ij} \neq
0$ only if $\om \lam_i = \lam_j$. In this case an easy computation
then yields that $B_{ij} = DX$, where $D$ is the $k_i \times k_i$
diagonal matrix $D =\diag(1,\om,\ldots,\om^{k_i-1})$ and $X$ is
a rectangular Toeplitz matrix (i.e. with equal elements on each
diagonal) such that all elements in the first column below position $(1,1)$
and all elements in the last row to the left of position
$(k_i,l_j)$ are $0$.

We note that in this way we have not obtained a canonical form for the pair
$(A,B)$ under simultaneous similarity, as in general there will be
similarities that leave $A$ invariant but change~$B$.

\begin{example}
Suppose $\alpha \neq 0$ and let
$$
A=\mat{ccccc} \lam & 1 & 0 & 0 & 0 \\
0 & \lam & 1 & 0 & 0 \\ 0 & 0 & \lam &0 & 0 \\
 0 & 0 & 0 & \alpha \lam & 1 \\ 0 & 0 & 0 & 0
& \alpha \lam  \rix.
$$
 If $\lam \neq 0$ and $\alpha \neq -1$ then $B$ is of the form
$$
B=\mat{ccccc} 0 & 0 & 0 & 0 & 0 \\
0 & 0 & 0 & 0 & 0 \\ 0 & 0 & 0 &0 & 0 \\
 0 & x_1 & x_2 & 0 & 0 \\ 0 & 0 & \alpha x_1 & 0 & 0
\rix.
$$
If $\lam \neq 0$ and $\alpha = -1$, then we obtain
$$
B=\mat{ccccc} 0 & 0 &
0 & y_1 & y_2  \\
0 & 0 & 0 & 0 & \alpha y_1 \\ 0 & 0 & 0 & 0 & 0 \\
 0 & x_1 & x_2 & 0 & 0 \\ 0 & 0 & \alpha x_1 & 0 & 0
\rix,
$$
 while, if $\lam  = 0$, then
$$
B=\mat{ccccc} u_1 & u_2 & u_3 & y_1 & y_2  \\
0 &  \alpha u_1 & \alpha u_2 & 0 &\alpha y_1 \\ 0 & 0 & \alpha^2 u_1 & 0 & 0 \\
 0 & x_1 & x_2 & v_1 & v_2 \\ 0 & 0 & \alpha x_1 & 0 & \alpha v_1
\rix.
$$
\end{example}

If $B$ is nonsingular and $A$ is not
nilpotent, then every row and every column of $B$ must contain at
least one nonzero element. Thus, if $\lam$ is a nonzero eigenvalue of $A$ so
is $\alpha \lam$. Since the number of eigenvalues is finite and $A$
has a nonzero eigenvalue, it follows that $\alpha$ is a root of
unity. Moreover, if  $J_i(\lam_i)$ is the
Jordan block of largest size in $A$, then using the fact that every
row and column of $B$ has at least one nonzero element, it follows that there
is a block for $\alpha \lam_i$ of equal size. Thus, we conclude that
the maximal size of a Jordan block in $A$ is the same for each
nonzero eigenvalue.

\begin{theorem} \label{jform}
Let $A$ and $B$ be nonsingular satisfying $AB = \alpha BA$, where $\alpha
\neq 1$ is non-zero. Then $\alpha $ is a primitive $p$-th root of $1$
for some $p$ and  the Jordan form of $A$ may be written as
$$ J = \diag(K_1, \ldots, K_q)$$
where each $K_k,\ k = 1,\ldots,q$, is the direct sum of Jordan
blocks of the same size belonging to $\lam, \alpha \lam, \ldots,
\alpha^{p-1} \lam$.

The same structure of the Jordan form holds for $B$.
\end{theorem}
\proof From the assumption of the theorem, $A$ is similar to $\omega A$.
Hence the number of Jordan blocks corresponding to any eigenvalue
$\lambda$ of $A$ and their sizes coincide with the number and sizes
of the Jordan blocks corresponding to the eigenvalue $\omega \lambda$.
The proof for $B$ follows by exchanging the roles of $A$ and $B$.
\eproof

Note further that, if $AB = \alpha BA$,
then $AB \cdot B = \alpha B \cdot AB$. Thus, if $A,B$ are
nonsingular quasi-commutative matrices, then $BA$ also satisfies
the conclusions of Theorem (\ref{jform}).

This finally brings us to the question of Wielandt on classification
of quasi-commutative pairs.
It was to a large extent answered by Drazin
already in 1951, although Wielandt was apparently unaware of his results.
In~\cite{Dra51}, Drazin obtained the following pre-normal form for
pairs of quasi-commutative matrices.
\begin{theorem}\label{drazin1}
If $A,B$ are $n\times n$ matrices satisfying an equation of the from
$AB=\omega BA$, then either
\begin{enumerate}
\item [(i)] $A,B$ can be simultaneously reduced to triangular form by a
similarity transformation,
or
\item [(ii)] there is an integer $r$ ($0\leq r \leq n-2$) such that
$A,B$
can be reduced, by the same similarity transformation, to the forms
\begin{equation}\label{drznf}
\mat{cc} S & X \\ 0 & A_r \rix, \qquad
\mat{cc} T & Y \\ 0 & B_r \rix,
\end{equation}
where $S,T$ are triangular $r\times r$ matrices, and $A_r,B_r$ are
nonsingular $(n-r)\times (n-r)$ matrices.
\end{enumerate}
\end{theorem}

Furthermore, Drazin also proves the following Theorem:
\begin{theorem}\label{drazin2}
In Theorem \ref{drazin1}, if (i) holds with $\omega\neq 1$, then
each of $AB$, $BA$ is nilpotent, and $A,B$ have between them at least
$n$ zero eigenvalues.
If, however, (i) is false, then $\omega$ is necessarily a
primitive root
of unity, and the order $k$ of $\omega$ must divide
$n-r$. Further, in this case, $ST$ and $TS$ are both nilpotent, and the
reduction of $A,B$ can be effected in such a way that $A_r$
can be reduced, by the same similarity transformation, to the forms
\begin{equation}\label{anrf}
\mat{cccc} a &&&\\ & \omega a && \\ &&\ddots&
 \\ &&& \omega^{k-1} a \rix,
\end{equation}
where $a$ is a non-singular square matrix of order $(n-r)/k$; then the
most general form of $B_r$ is
\begin{equation}\label{bnrf}
\mat{ccccc} 0& 0 &\ldots  &0& b_1\\
b_2 & 0 & \ldots &0 & 0 \\ &&\ddots&& \\
0 & 0 & \ldots & b_k & 0 \rix,
\end{equation}
where $b_1,\ldots b_k$ are arbitrary non-singular matrices of
order$(n-r)/k$;, subject to the relations
$b_i a=a b_i$, $(i=1,2,\ldots, k$)
are triangular $r\times r$ matrices, and $A_r,B_r$ are
nonsingular $(n-r)\times (n-r)$ matrices.
\end{theorem}

Drazin's formulas do not give a canonical form, however. Indeed, first of
all some further reduction of $A,B$ already in the form~(\ref{drznf}) is
possible. We have already seen that the equation $M_1X=XM_2$ has
only the trivial solution $X=0$, whenever the spectra of $M_1$ and $M_2$
do not intersect \cite{Gan59}.
This implies that we can decompose an arbitrary quasi-commutative pair
($A$, $B$) as
$$ A=\diag(\widetilde A_0,\widetilde A_1,\widetilde A_2,\ldots, 
\widetilde A_m), \qquad
B=\diag(\widetilde B_0,\widetilde B_1,\widetilde B_2,\ldots, 
\widetilde B_m),  $$
where the spectrum of $\widetilde A_0$ consists of zero only, and the spectra
of $\widetilde A_i$, $i=1, \ldots, m$
consist of distinct chains $\{ \lambda_i (\neq 0),
\omega \lambda_i, \ldots,  \omega^{k-1}\lambda_i \}$. (Note that
$\omega$ is a primitive $k$-th root of unity.) Then
each of the pairs $(\widetilde A_i,\widetilde B_i)$, $i=0, \ldots, m$
decomposes in the same way according to the spectrum of $\widetilde B_i$.
All together, we get a block-diagonalization
$$ A=\diag(A_0,A_1,A_2,\ldots, A_n), \qquad
B=\diag(B_0,B_1,B_2,\ldots, B_n)$$
such that each pair $(A_i,B_i)$ is of one of the following 4 types
according to the spectra $\sigma(A_i)$ and $\sigma(B_i)$:
$$ \begin{array}{llll}
 & \hbox{\rm Type I:} & \sigma(A_i)=\{0\}, & \sigma(B_i)=\{0\}, \\
 & \hbox{\rm Type II:} & \sigma(A_i)=\{ 0\}, &
\sigma(B_i)= \{\mu_i(\neq 0),  \omega \mu_i, \ldots, \omega^{k-1}\mu_i\}, \\
 & \hbox{\rm Type III:} & \sigma(A_i)=\{\lambda_i(\neq 0),
\omega \lambda_i, \ldots,
\omega^{k-1}\lambda_i\}, & \sigma(B_i)=\{0\}, \\
& \hbox{\rm Type IV:} & \sigma(A_i)=\{\lambda_i(\neq 0),
\omega \lambda_i, \ldots,
\omega^{k-1}\lambda_i\}, & \sigma(B_i)=\{\mu_i(\neq 0),  \omega \mu_i, \ldots,
\omega^{k-1}\mu_i\}.
\end{array} $$
Now, by Drazin's Theorem, each pair of type II can be put in the
form~((\ref{bnrf}),(\ref{anrf})) (notice
the order of matrices) and each pair of type III or IV to the
form~((\ref{anrf}),(\ref{bnrf})). Moreover, for a pair of type IV one can
assume that all submatrices $b_i$ in~(\ref{bnrf}), except for one
($b_1$ say), are equal to the identity. To achieve this, simply
use the transformation
$T=\diag(I,b_2^{-1},(b_2 b_3)^{-1}, \ldots, (b_2\cdots b_n)^{-1}).$
These form being fixed, the only further similarity transformations allowed
that do not destroy the identity blocks are
of the form $V=\diag (V_1, \ldots, V_1)$ with identical diagonal blocks
$V_1$ of the same size as the submatrix $a$. Therefore, the representation
problem for pairs of type IV reduces to the representation problem of
commuting matrix pairs $(a,b_1)$ under simultaneous similarity.
Conversely, the representation problem for commuting pairs of matrices under
simultaneous similarity reduces to the representation of quasi-commuting
pairs of type II, II or IV.  Indeed, suppose that
matrices $M$ and $N$ commute and are not both nilpotent.
By using the transformation to Jordan canonical form and
splitting the problem into subproblems, we may
 assume w.l.o.g. that $M$ has only one eigenvalue.
Moreover, since $M-\lambda I$  and $N -\mu I$
commute if and only if  $M$ and $N$ commute, we may assume at least one of
$M$ or $N$ to be nonsingular. Using $M$ as $a$ and and $N$ as $b_1$
in~(\ref{anrf}),~(\ref{bnrf}), and setting all other $b_i$'s to be $I$, we
obtain a quasi-commutative pair of type II, III or IV (depending on whether
$M$, $N$ or neither is chosen to be nilpotent). Since all transformations
preserving this form of the pair $(M,N)$ must look like $V=\diag
(V_1, \ldots, V_1)$, the problem of representing $(A,B)$ and that of
representing $(M,N)$ coincide. Note that this argument fails for pairs
$(M,N)$ where both matrices are nilpotent, since then it is no longer
true that $V$ has to be of the form $\diag(V_1, \ldots, V_1)$.

Our discussion can be summarized as follows.

\begin{theorem}
The problem of representation under simultaneous similarity for
quasi-commutative pairs is equivalent to the problem of representation
under simultaneous similarity for all commuting pairs. Moreover, the
latter is already equivalent to the problem of representation for
quasi-commutative pairs of type II, III or IV.
\end{theorem}

We do not know whether this result also holds for pairs of type I,
i.e., whether the problem of representation of quasi-commutative
pairs of type I is also equivalent to the problem of representation
of commuting pairs.


We now have an occasion to make a detour in the fascinating topic
of simultaneous similarity of commuting matrices.

M.~Gelfand and V.~A.~Ponomarev~\cite{GelP69} showed that the simultaneous
similarity problem of any n-tuple of matrices is equivalent to the simultaneous
similarity problem for a pair of commuting matrices. The seemingly hopeless
problem was later taken up by S.~Friedland, who showed in~\cite{Fri83} how
to find a finite number of invariants which will characterize a orbit of a
 pair $(A,B)$ under simultaneous similarity up to a finite ambiguity, which
means that these invariants may characterize a finite number of similarity
orbits. For a fixed dimension $d$, Friedland decomposes the variety of pairs
of square matrices in finitely many subsets locally closed under simultaneous
similarity (so subvarieties). For each of such subvarieties $Z$,
he gives a rational map $f$ from $Z$ into a finite-dimensional vector
space $V$, so that the pre-images under $f$ (of points in $V$) consist of
finitely many orbits of matrix pairs. With that, $f$ and $V$ depend strongly
on $Z$, while, for a fixed $f$, one can give an upper bound on the number of
conjugation classes in each pre-image.
Friedland's method was later refined by K.~Bongartz in~\cite{Bon95}.
He modified Friedland's construction (by changing $Z$, $f$ and $V$)
so that the pre-images under $f$ are exactly the individual orbits of
pairs of matrices.

In other words, given two pairs $(A,B)$ and $(C,D)$ of matrices, they are
simultaneously similar  to each other if and only if  they lie in the same
$Z$ and have the same image under $f$. This provides, at least in principle,
a complete answer to the problem of simultaneous similarity, i.e., a
a decision algorithm via rational computations, but no
readily available normal forms.

We should also mention that one of the abstract versions of
this problem is to find
all isomorphism classes of cyclic modules of finite length over the commutative
polynomial ring $R = \C[x,y]$. A pair of commuting $n\times n$-matrices $A$,
$B$ defines an $R$-module structure on $\C^n$ by letting $x$ and $y$ be
multiplication by $A$ and $B$, respectively.

\section{The converse to Potter's Theorem}\label{potterconv}
Having studied the decomposition of quasi-commutative matrices
into blocks, we now discuss
Wielandt's second question  whether the converse to Potter's theorem holds
for every irreducible block, i.e., whether the relation~(\ref{fermat})
where $s=t=1$ implies
\eq{qc} A B =\omega BA,\en
 where $\omega$ is the $q$-th root of unity.

If $q=2$, this strong version of the converse indeed holds.

\begin{proposition}\label{q=2}
A pair $(A,B)$ is quasi-commutative with $q=2$ if and only
if~(\ref{fermat}) holds with $s,t=1$.
\end{proposition}
\proof
The condition $A^2+B^2=(A+B)^2$ is equivalent to $AB=-BA$.
\eproof

However, the converse is in general not true, even if~(\ref{fermat}) is
assumed to hold for all values of $s$ and $t$, as the following example shows.

\begin{example}\label{cex1}\rm

It is in general not true that if~(\ref{fermat}) holds for some $q$ and all
$s$, $t$, then~(\ref{qc}) holds with $\omega$ the $q$-th root of unity.

Consider the case $n=3$, $q=3$ and let  $\lambda\in \C$
be not equal to $0$, $-1$ or one of the two primitive
$3$rd roots of unity. For the pair of matrices
\eq{ceconv}
A=\mat{ccc} 0 & 0 &\frac 1\lambda \\
1 & 0 & 0 \\ 0 & -\frac {\lambda+1}{\lambda} & 0
\rix,\quad B=\mat{ccc} 0 & 1 & 0 \\ 0 & 0 & 1 \\ 1 & 0 & 0 \rix.
\en
we have
\eq{ecom} BA=EAB\en
 where
$$
E=\mat{ccc} \lambda & 0 & 0 \\ 0 & -\frac {\lambda+1}{\lambda} & 0 \\
0 & 0 & -\frac 1 {\lambda+1} \rix.
$$
Moreover, since $E$ is invertible, it follows that
\eq{newcomm} A(E^{-1} +I)=-EA, \quad (E^{-1}+I)B=-BE.
\en
But (\ref{ecom}) and (\ref{newcomm}) imply that
$$
(A+tB)^3= A^3+t^3 B^3
$$
for all $t$. However, since $E$ has $3$ distinct eigenvalues it follows
that~(\ref{qc}) does not hold. Also, if the pair $(A,B)$ is replaced by
$(\widetilde A, \widetilde B)\eqbd T(A,B)T^{-1}$, then
$$ \widetilde B \widetilde A=\widetilde E \widetilde A \widetilde B \qquad
{\rm where} \quad \widetilde E \eqbd TET^{-1}.$$
Since $E$ and $\widetilde E$ have the same spectrum,~(\ref{qc})
does not hold for the pair $(\widetilde A, \widetilde B)$ either.
In other words, the pair $(A,B)$ cannot be reduced to a direct sum
of quasi-commutative pairs.

Note that in this example both matrices $A$ and $B$ are nonsingular.
\end{example}

If we assume that $s=t=1$ in~(\ref{fermat}), then we can produce
even $2\times 2$-counterexamples, e.g., with $q=3$.

\begin{example}\rm \label{ex4}
The  pair of matrices
\eq{2x2ce1}
A=\mat{cc} 1 & 1 \\ 0 & 2\rix,\qquad B=\mat{cc} -1 & 0 \\ 7 & -2\rix
\en
satisfies (\ref{fermat}) with $s,t=1, q=3$ but, since
$AB=\mat{cc} 6 & -2 \\ 14 & -4\rix$ and $BA= \mat{cc} -1 & -1 \\ 7 &
-4 \rix$,
 the pair is not quasi-commutative.

The  pair of matrices \eq{2x2ce} A=\mat{cc} 0 & 1 \\ 0 &
2\rix,\qquad B=\mat{cc} -2 & 2 \\ 0 & 0\rix \en is even triangular
and (\ref{fermat}) holds with $s,t=1, q=3$ but, since $AB=0$ and
$BA= \mat{cc} 0 & 2 \\ 0 & 0 \rix$, the pair is not
quasi-commutative.

One can check that the pairs in both of these examples cannot be
decomposed into direct sums of quasi-commutative matrices either.
\end{example}

Drazin's pre-normal form for quasi-commutative matrices also suggests
the question whether the converse to Potter's theorem holds at
least for the pairs of matrices of the form~(\ref{drznf})(i). The
following examples  demonstrate that it is not so.
\begin{example}\rm \label{ex5}
Let $n=q=3$, let $\omega$ be either of the two primitive 3-rd roots
of unity and let $x_i$,
$i=1, 3$ be arbitrary nonzero numbers. Consider the triangular matrices
\eq{cetriangle}
A=\mat{ccc}1 & 0 & 0 \\
0 & \omega & 0  \\ 0 & 0 & \omega^2
\rix,\quad B=\mat{ccc} 0 & x_1 & x_2 \\ 0 & 0 & x_3
\\ 0 & 0 & 0 \rix. \en
Then $A^3=I$, $B^3=0$, and $(A+tB)^3=I$ for any $t$.
But
\eq{compare}
AB=\mat{ccc} 0 & \omega x_1 & \omega^2 x_2 \\
0 & 0 &  \omega^2 x_3  \\ 0 & 0 & 0
\rix,\quad BA=\mat{ccc} 0 & x_1 & x_2 \\ 0 & 0 & \omega x_3 \\
0 & 0 & 0 \rix,
\en
so the products $AB$ and $BA$ are not scalar multiples of each other.
Since the matrix $B$ is similar to a single Jordan block of
size $3$ corresponding to the eigenvalue $0$, the pair $(A,B)$
is also seen to be irreducible.
\end{example}

On the other hand, if a pair of block $r\times r$-matrices of the
form~(\ref{anrf})-(\ref{bnrf}) with
some $\omega$ satisfies~(\ref{fermat}), then necessarily $r=q$,
$\omega^q=1$,
and~(\ref{qc}) holds. Indeed, suppose a pair $(A,B)$ is block $r\times r$
in the form~(\ref{anrf})-(\ref{bnrf}). Then, by direct calculation, it
satisfies~(\ref{qc}), and, comparing the determinants on both sides, we
obtain $\omega^r=1$ . So, we only need to establish
that $q=r$. The relations~(\ref{fermat})
and~(\ref{qc}) together imply that $\omega^q=1$, hence $q=rm$ for some
natural number $m$. So, the pair
 $(\widetilde{A}\eqbd A^r,\widetilde{B} \eqbd B^r)$ satisfies the relation
\eq{newrel}\widetilde{A}^m+\widetilde{B}^m=(\widetilde{A}^m+\widetilde{B}^m).
\en
But the matrices $\widetilde{A}$ and $\widetilde{B}$ commute,
hence~(\ref{newrel})
implies that $m=1$. Thus, $q=r$.

However, the commutativity of the blocks $a$ and $b_i$
in~(\ref{anrf})-(\ref{bnrf}) does not follow automatically from the
relation~(\ref{fermat}), so here again the converse to the Potter's
result fails, as we show next.

\begin{example}\rm \label{ex6}
 Let $n=6$, $q=3$ and let
\eq{ceblocks}
A=\mat{ccc} A_1 & 0 & 0 \\
0 & \omega A_1 & 0  \\ 0 & 0 & \omega^2 A_1
\rix,\quad B=\mat{ccc} 0 & B_1 & 0 \\ 0 & 0 & B_1 \\ B_1 & 0 & 0 \rix,
\en
where
$$ A_1=\mat{cc} 1 & 1+{\omega \over 1-\omega^2}  \\
0 & \omega^2 \rix,\quad B=\mat{cc} 1 & 1  \\ 0 & \omega^2 \rix.
$$
Then $(sA+tB)^3=(sA)^3+(tB)^3$ for all $s$, $t$, but $B_1A_1=JA_1B_1$,
with  $$ J=\mat{cc} 1 & 1 \\
0 & 1 \rix,
$$
so, in particular, $A_1$ and $B_1$ do not commute and hence the
matrices $A$, $B$ are not quasi-commutative either.
\end{example}

In view of these counterexamples, it seems natural to pose
the  more general problem to characterize all classes of matrices for
which the equivalence
of (\ref{fermat}) and (\ref{qc}) holds.

\section{Identities satisfied by quasi-commutative matrices}\label{identities}

The first question Wielandt asked was which identities are satisfied
by quasi-commutative matrices. We now show that all
polynomial identities $f(x,y)=0$ that hold for all
quasi-commutative matrices belong to the ideal in $\C[x,y]$
generated by the polynomial $xy-\omega yx$.

\begin{theorem}
Let $\C[x,y]$ denote the ring of polynomials in non-commuting
indeterminates $x$, $y$
over the field $\C$ and let $\cal I$ denote the ideal of  $\C[x,y]$
generated by
the polynomial $g(x,y)=xy-\omega yx$ with $\omega^q=1$. Then
\eq{ideal} f(x,y)\in {\cal I} \en
if and only if the condition~(\ref{qc}) implies $f(A,B)=0$ in $\C^{q\times q}$.
\end{theorem}

\proof  One direction is obvious: any polynomial $f(x,y)\in {\cal
I}$ satisfies $f(A,B)=0$ for all quasi-commutative matrices $A$,
$B$.

To show the converse, first recall that the condition
$\omega^q=1$ implies that
there exist a pair of {\em nonsingular\/} matrices $A$,
$B\in \C^{q \times q}$
satisfying~(\ref{qc}).
Since the pair $(sA,tB)$ also satisfies~(\ref{qc}) for any
scalars $s$, $t$, we get
$$ f(sA, tB)=0.$$
Now interchange $A$ and $B$ using the relation~(\ref{qc})
as many times as to obtain a polynomial in the form
$$ f_1(sA, tB) \eqbd \sum_{i,j}c_{i,j} s^i t^j A^i B^j. $$
The polynomials $f(x,y)$ and $f_1(x,y)$ differ by some element of $\cal I$.
Now, since $f_1(sA, tB)=0$ and $s$ $t$ are independent scalars, each
term $c_{i,j} s^i t^j A^i B^j$ in the sum must equal zero. But as both
$A$ and $B$ are nonsingular, this shows that $c_{i,j}=0$. Thus, $f_1(x,y)$
is the zero polynomial and hence $f(x,y)\in {\cal I}$.
\eproof

\section{Further historical comments}\label{final}

We now address the last question asked by Wielandt, namely on the
work preceding that of H. S. A. Potter. This question turns out to
be also briefly answered by M. P. Drazin in~\cite{Dra51}.
Specifically, Drazin cites Cayley's paper~\cite{Cay58} where the
case $\omega=-1$ was considered and the works of
F.~Cecioni~\cite{Cec31}, S.~Cherubino~\cite{Che38} and
T.~Kurosaki~\cite{Kur41} devoted to the general case. [Biographies of
the two Italian mathematicians can be found at~\cite{ItWeb}.]
Cecioni's paper is a memoir summarizing and extending results on
quasi-commutative matrices known at that time. He proves a
condition on a matrix $A$ necessary and sufficient for the
equation $AX=\omega XA$ to have a nonzero solution $X$, describes
the structure of an arbitrary solution similarly to
Turnbull-Aitken~\cite[p.148]{TurA32}, and stops one step before
arriving at the formulas~(\ref{anrf})-(\ref{bnrf}) for a
quasi-commutative pair $(A,B)$ with $AB$ nonsingular. A slightly
different pre-normal form is derived by Cherubino~\cite{Che38}; he
also describes the structure of the algebra of matrices commuting
with a given matrix. The pair~(\ref{anrf})-(\ref{bnrf}) appears
also in Kurosaki~\cite{Kur41}, even in the reduced form (with all
$b_j$'s except for one equal to the identity), although not in a
formal statement. Kurosaki's main result (\cite[Theorem 4]{Kur41})
is a description of the group of all nonsingular matrices $P$
satisfying the equation $AP=cPA$ for some $c$ (depending on $P$)
and a fixed nonsingular matrix $A$.  Drazin in~\cite{Dra51} is
apparently more interested in simultaneous triangularization of a
quasi-commutative pair, hence obtains, in his remarkably short
paper, yet another reduced form.

\section*{Acknowledgements} We are indebted to H. S. A.~Potter and
M. P.~Drazin for further references, to S.~Friedland,
L.~Levy, and K.~Bongartz for discussions of the simultaneous similarity
problem, to R.~Askey for information on the $q$-binomial formulae,
to M.~Lorenz for drawing our attention to the connection to quantum planes,
and to M.~Benzi and N.~Laghi for their help with the works of Cecioni
and Cherubino.

\bibliographystyle{plain}

\end{document}